\documentclass[10pt]{amsart} 
\usepackage{amssymb,latexsym}

\newtheorem{theorem}{Theorem}

\newtheorem{corollary}[theorem]{Corollary}
\newtheorem{definition}[theorem]{Definition}
\newtheorem{conjecture}[theorem]{Conjecture}

\newtheorem{property}{Property}

\def\wnot{w_\mathrm{o}}

\def\pa{\partial}
\def\S{\mathfrak{S}}
\def\Q{\mathfrak{S}^q}

\def\H{{\rm H}}
\def\QH{{\rm QH}}

\def\l{\ell}
\def\E{\mathcal{ E}}

\def\s{\sigma}
\def\<{\langle\kern-.08cm\langle}
\def\>{\rangle\kern-.08cm\rangle}

\def\dunkl{\theta}

\def\Enp{\E_n^+}
\def\Eqn{{{\E}^q_n}}

\title[Quantum cohomology of the flag manifold]{Lecture notes on
quantum cohomology\\ of the flag manifold}

\date{January 15, 1999}

\author{Sergey Fomin}

\address{
Department of Mathematics,
Massachusetts Institute of Technology,
Cambridge, MA 02139
}

\thanks{Partially supported by the NSF (DMS-9700927).} 

\begin{document}

{\ }

\vspace{-1in}

\setcounter{property}{0}


\maketitle

This is an exposition of some recent developments related to the object in the title,
particularly the computation of the  Gromov-Witten invariants 
of the flag manifold~\cite{FGP} and the quadratic algebra approach~\cite{FK}. 
The notes are largely based on the papers~\cite{FGP} and~\cite{FK},
authored jointly with S.~Gelfand, A.~N.~Kirillov, and A.~Postnikov. 
This is by no means an exhaustive survey of the subject, but rather a casual
introduction to its combinatorial aspects. 

\section{Classical theory}
\label{sec:intro}

Let us briefly review the standard facts from the Schubert calculus of
the flag manifold; see \cite{fulton-yt} for details. 
Let $Fl_n$ be the variety of complete flags in~$\mathbb{C}^n$.
The cohomology ring $\H^*(Fl_n,\mathbb{Z})$  
can be described in two different ways. 
The first description, due to Borel~\cite{borel},
represents it as a quotient of a polynomial ring:
\begin{equation}
\label{eq:borel}
\H^*(Fl_n,\mathbb{Z})\cong \mathbb{Z}[x_1,\dots,x_n]/I_n,
\end{equation}
where $x_1,\dots,x_n\in \H^2(Fl_n,\mathbb{Z})$ are the first Chern classes of
$n$ standard line bundles on $Fl_n\,$, and $I_n$ is the ideal generated by
symmetric polynomials in $x_1,\dots,x_n$ without constant
term\footnote{This result, as well as several others below, extends to
the more general setup of the homogeneous space $G/B$ for a complex
semisimple Lie group~$G$.  
In these notes, we only treat the type~$A$ case, with~$G=SL_n\,$.}. 

The second description is based on the decomposition of $Fl_n$ 
into Schubert cells, indexed by the elements 
of the symmetric group~$S_n\,$. 
The corresponding cohomology
classes $\s_w\,$, $w\in S_n$ (the Schubert classes)
form an additive basis in $\H^*(Fl_n,\mathbb{Z})$. 

The elements of the quotient ring $\mathbb{Z}[x_1,\dots,x_n]/I_n$
which correspond to the Schubert classes under the isomorphism 
(\ref{eq:borel}) were identified by Bernstein, Gelfand, and
Gelfand~\cite{BGG} and Demazure~\cite{Dem}. 
Then Lascoux and Sch\"utzenberger~\cite{ls} introduced remarkable
polynomial representatives of the Schubert classes called Schubert polynomials.
These polynomials $\S_w\,$, $w\in S_n\,$, 
are defined as follows.

Let $s_i$ denote the adjacent transposition $(i~i+1)$.
For $w\in S_n\,$, an expression
$w=s_{i_1}s_{i_2}\dots s_{i_l}$ of minimal possible length is called a
\emph{reduced decomposition}. 
The number $l=\l(w)$ is the \emph{length} of~$w$. 
The symmetric group~$S_n$ acts on 
$\mathbb{Z}[x_1,\dots,x_n]$ by 
$w\, f=f(x_{w^{-1}(1)},\dots,x_{w^{-1}(n)})$. 
The \emph{divided difference operator} $\pa_i$ is defined by
$\pa_i\, f=(x_i-x_{i+1})^{-1}(1-s_i)\, f$. 
For any permutation $w$, the operator $\pa_w$ 
is defined by $\pa_w=\pa_{i_1}\pa_{i_2}\dots \pa_{i_l}$,
where $s_{i_1}s_{i_2}\dots s_{i_l}$ is a reduced decomposition for~$w$. 

Let $\delta=\delta_n=(n-1,n-2,\dots,1,0)$ and 
$x^\delta=x_1^{n-1}x_2^{n-2}\dots x_1$. For 
$w\in S_n\,$, the \emph{Schubert polynomial} $\S_w$ 
is defined by
$\S_w=\pa_{w^{-1}\wnot}x^\delta$, 
where $\wnot$ is the longest element in~$S_n\,$.
Equivalently, $\S_{\wnot}=x^\delta$, and 
$\S_{ws_i}=\pa_i\S_w$ whenever $\l(ws_i)=\l(w)-1\,$.
The following result is immediate from~\cite{BGG}.

\begin{theorem}
\label{th:BGG}
The Schubert polynomials represent Schubert classes under Borel's
isomorphism~{\rm(\ref{eq:borel})}. 
\end{theorem}

\section{Quantum cohomology}

The (small) \emph{quantum cohomology ring} $\QH^*(X,\mathbb{Z})$ of a smooth
algebraic variety~$X$ is a certain deformation 
of the classical cohomology; see, e.g., \cite{fulton-rahul} for
references and definitions. 
The additive structure of this ring is usually rather simple. 
For example, $\QH^*(Fl_n,\mathbb{Z})$  is canonically isomorphic, 
as an abelian group, to the tensor product 
$\H^*(Fl_n,\mathbb{Z})\otimes \mathbb{Z}[q_1,\dots,q_{n-1}]$, where
the $q_i$ are formal variables (deformation parameters). 
The multiplicative structure of the quantum cohomology is however
deformed comparing to $\H^*(Fl_n,\mathbb{Z})$, 
and specializes to it in the classical limit $q_1=\cdots=q_{n-1}=0$.
The multiplication in $\QH^*(Fl_n,\mathbb{Z})$ is given by 
\begin{equation}
\label{eq:2.3.1}
\s_u * \s_v = \sum_w \sum_{d=(d_1,\dots,d_{n-1})}q^d
\left<\s_u,\s_v,\s_w\right>_d \,\s_{\wnot w}\, , 
\end{equation}
where the $\left<\s_u,\s_v,\s_w\right>_d$ are the (3-point, genus~$0$)
\emph{Gromov-Witten invariants} of the flag manifold,
and $q^d=q_1^{d_1}\cdots q_{n-1}^{d_{n-1}}$.  
Informally, these invariants count equivalence class\-es 
of rational curves in $Fl_n$
which have multidegree $d=(d_1,\dots,d_{n-1})$ 
and pass through given Schubert varieties.
In order for an invariant to be nonzero,  the condition
$
\l(u)+\l(v)+\l(w)=\binom{n}{2}+2\sum_{i=1}^{n-1}d_i 
$
has to be satisfied. 
The operation~$*$ defined by (\ref{eq:2.3.1}) is associative
\cite{li-tian, RT}, and obviously commutative.

The quantum analog of Borel's theorem was 
obtained by Givental and Kim \cite{giv-kim, Kim1, Kim2, Kim3} and 
Ciocan-Fontanine~\cite{ciocan} who showed that
\begin{equation}
\label{eq:H*q}
\QH^*(Fl_n,\mathbb{Z})\cong
P_n/I^q_n\ ,
\end{equation}
where $P_n=\mathbb{Z}[q_1,\dots,q_{n-1}][x_1,\dots,x_n]$, 
the $x_i$ are the same as before, 
and $I^q_n$ is the ideal generated by 
the coefficients~$E_1^n,\dots,E_n^n$ 
of the characteristic polynomial 
\begin{equation}
\label{eq:E_n}
\det(1+\lambda G_n) = \sum_{i=0}^n E_i^n \lambda^i 
\end{equation}
of the matrix
\begin{equation}
\label{eq:G_n}
G_n=
\left(
\begin{array}{ccccc}
  x_1    & q_1    & 0      & \cdots  & 0      \\
  -1     & x_2    & q_2    & \cdots  & 0      \\
  0      & -1     & x_3    & \cdots  & 0      \\
  \vdots & \vdots & \vdots & \ddots  & \vdots \\
  0      & 0      & 0      & \cdots  & x_n
\end{array}
\right)\ .
\end{equation}
(These coefficients are called \emph{quantum elementary
symmetric functions}.)
More precisely, let us identify the polynomial $x_1+\dots+x_i\,$
with the Schubert class $\s_{s_i}\,$.
The quantum cohomology ring is then generated by the elements~$x_i\,$,
subject to the relations in the ideal~$I^q_n\,$.

\section{Quantum Schubert polynomials}

The above description of $\QH^*(Fl_n,\mathbb{Z})$ 
does not tell which elements on the right-hand side of (\ref{eq:H*q})
correspond to the Schubert classes. 
The main goal of \cite{FGP} was to give the quantum analogues 
of the Bernstein--Gelfand--Gelfand theorem and
the Schubert polynomials construction of Lascoux and Sch\"utzenberger.
This allowed us to design algorithms for computing
the Gromov-Witten invariants for the flag manifold. 
Our approach relied on some of the most basic
properties of the quantum cohomology, 
which can be expressed in elementary terms (see below).

Let~$A_n$ denote the vector space spanned
by the classical Schubert polynomials.
Another basis of $A_n$ is formed by the monomials 
$x_1^{a_1}x_2^{a_2}\dots x_{n-1}^{a_{n-1}}$ dividing the staircase
monomial~$x^\delta$.
The space~$A_n$ is complementary to the ideal~$I_n\,$,
and also to the quantized ideal~$I^q_n\,$. 

The \emph{quantum Schubert polynomial} $\Q_w$ is defined as the unique
polynomial in $A_n$ that belongs to the coset  
modulo~$I^q_n$ representing the Schubert class~$\s_w$ 
under the canonical isomorphism (\ref{eq:H*q}). 
The primary goal of~\cite{FGP} was to algebraically identify these
polynomials.  

\section{Axiomatic characterization} 

The following properties of the quantum Schubert polynomials are
directly implied by their definition.

\begin{property}
\label{property:homogeneous}
$\Q_w$ is homogeneous of degree $\l(w)$,
assuming $\deg(x_i)\!=\!1$, $\deg(q_j)\!=\!2$.
\end{property}

\begin{property}
\label{property:specialization}
Specializing $q_1=\dots=q_{n-1}=0$ yields $\Q_w=\S_w\,$.
\end{property}

\begin{property}
\label{property:staircase}
$\Q_w\,$ belongs to the span~$A_n$ of the classical Schubert polynomials.
\end{property}

It follows that the~$\Q_w$ form a linear basis in~$A_n\,$,
and that the transition matrices between the bases
$\{\Q_w\}$ and $\{\S_w\}$
are unipotent triangular, with respect to any linear ordering
consistent with~$\l(w)$.

The next property reflects the fact that the 
Gromov-Witten invariants of the flag manifold are nonnegative integers.

\begin{property}
\label{property:structure-constants}
Consider any product of polynomials~$\Q_w\,$.
Expand it (modulo $I^q_n$)
in the linear basis~$\{\Q_w\}$. 
Then all coefficients in this expansion are polynomials in the $q_j$
with nonnegative integer coefficients.
\end{property}


The following result is a restatement
of formula~(3) in~\cite{ciocan}. 

\begin{property}
\label{property:E-ciocan}
For a cycle $w=s_{k-i+1}\cdots s_{k}\,$, 
we have $\Q_w=E^k_i\,$.
\end{property}


\begin{theorem}
\label{th:char} \cite{FGP} 
The polynomials $\Q_w$ are uniquely determined by 
Properties~\ref{property:homogeneous}--\ref{property:E-ciocan}.
\end{theorem}

We conjecture in~\cite{FGP} that Property~\ref{property:E-ciocan},
which is the only property stated above that does not trivially follow
from the quantum-cohomology definition of the $\Q_w\,$,
is not actually needed 
to uniquely determine the quantum Schubert polynomials.

The next two sections provide constructive descriptions of these polynomials.




\section{Quantum polynomial ring}
\label{sec:q-mult}

For $k=1,2,\dots$, define the operator $X_k$ acting
in the polynomial ring by
\begin{equation}
\label{eq:Qxk}
X_k=x_k\,-\,\sum_{i<k}q_{ik}\pa_{(ik)}+\sum_{j>k}q_{kj}\pa_{(kj)}\ ,
\end{equation}
where 
$\pa_{(ij)}$
is the divided difference operator which 
corresponds to the transposition $t_{ij}\,$, and $q_{ij}=q_iq_{i+1}\dots q_{j-1}$.
(We will always assume~$i<j$.)

\begin{theorem}
\label{th:commute}
\cite{FGP}
The operators $X_i$ commute pairwise, and generate a free commutative ring.
For any polynomial $g\in P_n\,$,
there exists a unique operator $G\in \mathbb{Z}[q_1,\dots,q_{n-1}][X_1,\dots,X_n]$
satisfying~$g=G(1)$.
\end{theorem}

\noindent
(Here $G(1)$ denotes the result of applying $G$ to the polynomial~$1$.)

For a polynomial $g\in P_n\,$, the polynomial $G$ given by $g=G(1)$ is
called the \emph{quantization} of~$g$. 
The bijective correspondence $g\leftrightarrow G$ 
between $P_n$ and $\mathbb{Z}[q_1,\dots,q_{n-1}][X_1,\dots,X_n]$ 
is by no means a ring homomorphism.
Identifying the two spaces via this bijection,
we obtain an alternative ring structure on~$P_n\,$. 
The multiplication thus defined is called 
\emph{quantum multiplication} and denoted by~$*$; 
it coincides with the usual multiplication  
in the classical limit.

Recall that $I_n\subset P_n$ is the ideal 
generated by the elementary symmetric functions 
$e_i=e_i(x_1,\dots,x_n)$, $i=1,\dots,n$.
It can be checked that $I_n$ is also an ideal 
with respect to the quantum multiplication 
(i.e., $I_n$ is an invariant space for the operators
$X_1,\dots,X_n$ acting in~$P_n$). 

We are now going to relate our quantum multiplication
to the quantum cohomology of the flag manifold.
First we verify that for $i=1,\dots,n$, the quantization of the elementary
symmetric function $e_i(x_1,\dots,x_n)$, 
is the quantum elementary symmetric function $E^n_i$ defined by~{\rm 
(\ref{eq:E_n})}. 
As a corollary, the quantization map bijectively maps the ideal $I_n$
onto the Given\-tal-Kim ideal~$I^q_n\,$. 
Thus the quotient $P_n/I_n\,$, with the quantum multiplication~$*$
defined above, is canonically isomorphic
to the quotient ring $P_n/I_n^q$ 
(hence to $\QH^*(Fl_n,\mathbb{Z})$).
In fact, more is true.

\begin{theorem}
\label{th:main}
\cite{FGP}
The canonical isomorphism between the quotient space $P_n/I_n$
and the classical cohomology of the flag manifold
is also a ring isomorphism between $P_n/I_n\,$,
endowed with quantum multiplication defined above in this section,     
and the quantum cohomology ring of the flag manifold.
\end{theorem}

In other words, the identification of
the (classical) Schubert polynomials with the corresponding
Schubert classes translates the quantum multiplication
defined in this section into the multiplication
in the quantum cohomology ring.

The quantum Schubert polynomial
$\S_w^q$ is the quantization
of the ordinary Schubert polynomial $\S_w\,$,
in the sense of the above construction.
In other words, $\S_w^q$ is uniquely determined by
$\S_w^q(X_1,X_2,\dots)(1)=\S_w(x_1,x_2,\dots)$. 
It follows that the quantum multiplication
of ordinary Schubert polynomials translates into the ordinary
multiplication of the corresponding quantum Schubert polynomials.

\section{Standard monomials}
\label{sec:stand-mon}

Let $e_i^k$ denote the elementary symmetric function of degree~$i$ 
in the variables $x_1,\dots,x_k\,$.
The \emph{standard elementary monomials} are defined by the formula 
\begin{equation}
\label{eq:standard} 
e_{i_1\dots i_{n-1}} = e_{i_1}^1 \dots e_{i_{n-1}}^{n-1}  \ ,
\end{equation}
where we assume $0\leq i_k\leq k$ for all $k$. 
It is well known (and easy to prove) that the
polynomials~(\ref{eq:standard}), for a fixed~$n$,  
form a linear basis in the space~$A_n$ spanned by the Schubert polynomials
for~$Fl_n\,$. 
Each Schubert polynomial $\S_w$ is thus uniquely expressed as a
linear combination of such monomials.

Let $G_k$ denote the $k$th leading principal minor of the matrix $G_n$
given by~(\ref{eq:G_n}). 
The \emph{quantum standard elementary monomial} is defined by
\begin{equation}
\label{eq:Em}
E_{i_1\dots i_{n-1}} = E_{i_1}^1 \dots E_{i_{n-1}}^{n-1} \ ,
\end{equation}
where $E^k_i=E_i(X_1,\dots,X_k)$ denotes the coefficient  of $\lambda^i$
in the characteristic polynomial $\chi(\lambda)=\det(1+\lambda G_k)$
of~$G_k\,$. 

\begin{theorem}
\label{th:intro}
\cite{FGP}
The quantum Schubert polynomial 
$\S^q_w$ is obtained by
replacing each standard monomial {\rm(\ref{eq:standard})} in the expansion of~$\S_w$ 
by its quantum analogue~{\rm (\ref{eq:Em})}.
\end{theorem}

The expansions of Schubert polynomials in terms of the standard
monomials can be computed recursively top-down in the weak
order of~$S_n\,$, starting from $\S_{\wnot}=e_{12\dots {n-1}}\,$.
Namely, use the basic divided difference recurrence for the $\S_w$ 
together with the rule for computing a divided difference
of an elementary symmetric function,
the Leibnitz formula for the~$\pa_i\,$,
and the corresponding straightening procedure. 
Having obtained such an expansion for $\S_w\,$, ``quantize'' each term in it to
obtain~$\S_w^q\,$.
In the special case $n=3$, this produces results shown in Figure~\ref{fig:2}. 

\setlength{\unitlength}{0.5mm}
\thicklines
\begin{figure}[ht]
\begin{center}  
\begin{picture}(40,70)(0,-5)
    \put(20,0){\circle*{3}} \put(0,40){\circle*{3}}
    \put(0,20){\circle*{3}} \put(20,60){\circle*{3}}
    \put(40,20){\circle*{3}} \put(40,40){\circle*{3}}
    \put(20.4,-1){\line(1,1){20}}
    \put(19.6,1){\line(1,1){20}}
    \put(20,0){\line(-1,1){20}}
    \put(-1,20){\line(0,1){20}}
    \put(1,20){\line(0,1){20}}
    \put(40,20){\line(0,1){20}}
    \put(0,40){\line(1,1){20}}
    \put(39.4,39.4){\line(-1,1){20}}
    \put(40.6,40.6){\line(-1,1){20}}
    \put(17,-9){$1$}
    \put(-10,19){$s_1$}
    \put(45,19){$s_2$}
    \put(-19,39){$s_1 s_2$}
    \put(45,39){$s_2s_1$}
    \put(15,67){$\wnot$}
  \end{picture}\qquad\qquad\qquad\qquad\qquad
\begin{picture}(40,70)(0,-5) 


\put(20,0){\circle*{3}}
\put(0,40){\circle*{3}}
\put(0,20){\circle*{3}}
\put(20,60){\circle*{3}}
\put(40,20){\circle*{3}}
\put(40,40){\circle*{3}}


\put(20.4,-1){\line(1,1){20}}
\put(19.6,1){\line(1,1){20}}
\put(20,0){\line(-1,1){20}}
\put(-1,20){\line(0,1){20}}
\put(1,20){\line(0,1){20}}
\put(40,20){\line(0,1){20}}
\put(0,40){\line(1,1){20}}
\put(39.4,39.4){\line(-1,1){20}}
\put(40.6,40.6){\line(-1,1){20}}


\put(17,-9){$1$}
\put(-12,20){$x_1$}
\put(43,20){$x_1+x_2$}
\put(-33,40){$x_1x_2+q_1$}
\put(43,40){$x_1^2-q_1$}
\put(4,65){$x_1^2x_2+q_1x_1$}

\end{picture}
\end{center} 
\caption{Quantum Schubert polynomials for $S_3$}
\label{fig:2}
\end{figure}
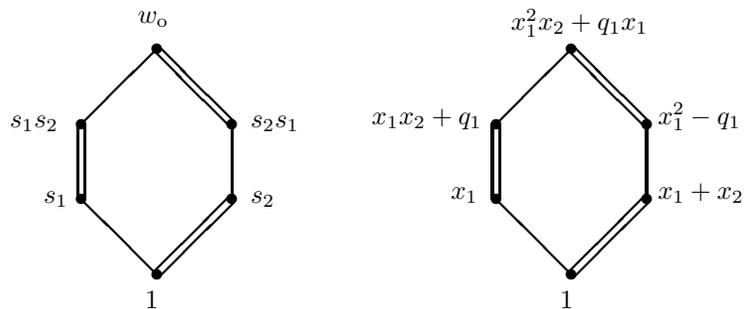

\section{Computation of the Gromov-Witten invariants}

The space~$A_n$ spanned by the Schubert polynomials for $S_n$ 
can be described as the set of normal forms
for the ideal~$I^q_n\,$, with respect to certain term order.
This allows one to use Gr\"obner basis techniques (see, e.g.,
\cite{sturmfels}) to construct efficient algorithms for computing 
the Gromov-Witten invariants of the flag manifold.

\begin{definition} 
\label{def:can-form}
{\rm Let us choose the \emph{total degree -- inverse lexicographic 
term order} on the monomials $x_1^{a_1}\cdots x_n^{a_n}\,$.
In other words, we first order all monomials by the total degree
$\sum_i a_i\,$,
and then break the ties by using the inverse lexicographic order 
$x_1<x_2<x_3<\dots\,$.
This allows us to introduce the \emph{normal},
or fully reduced \emph{form} of any polynomial with respect 
to the ideal $I^q_n$ and the term order specified above.
This normal form can be found, e.g., via the Buchberger algorithm
employing the corresponding Gr\"obner basis of~$I^q_n\,$.
}
\end{definition}

\begin{theorem} 
\label{th:can-form}
\cite{FGP} 
Choose a term order as in Definition~\ref{def:can-form}. 
Then the reduced minimal Gr\"obner basis for the ideal $I^q_n$
consists of the polynomials 
$\det{\left(E_{j-i+1}^{n-i+1}\right)_{i,\,j=1}^k}\,$, 
for $k=1,\dots,n$.
The normal form of any polynomial $F\in P_n\,$,
lies in the space~$A_n\,$.
\end{theorem} 

For a polynomial $F\in P_n\,$, let 
\begin{equation*}
\label{eq:top-coeff}
\< F\>
= {\rm ~coefficient~of~} x^\delta {\rm ~in~the~normal~form~of~} F\,.
\end{equation*}
Equivalently, $\< F \>$ is the coefficient of $\S^q_{\wnot}$
in the expansion of $F$ (modulo~$I^q_n$) in the basis of quantum
Schubert polynomials, since $\S^q_{\wnot}$
is the only basis element that involves
the staircase monomial~$x^\delta$.
The definition (\ref{eq:2.3.1}) implies that 
\[
\< \S^q_{w_1}\cdots \S^q_{w_k}\> 
= \sum_d q^d \langle \s_{w_1},\dots,\s_{w_k} \rangle_d\ ,
\]
the generating function for the Gromov-Witten invariants.
We thus arrived at the following result.

\begin{theorem}
\label{th:g-w}
\cite{FGP} 
A Gromov-Witten invariant 
$\langle\s_{w_1}\,,\dots,\s_{w_k}\rangle_d$ of the flag manifold 
is the coefficient of the monomial
$q^d\,x^\delta$ in the normal form 
(in the sense of Definition~\ref{def:can-form})
of the product of quantum Schubert polynomials
$\S^q_{w_1}\cdots \S^q_{w_k}\,$.
\end{theorem} 

\section{Quadratic algebras}

Another approach to the study of the cohomology ring---ordinary or
quantum---of the flag manifold was suggested in~\cite{FK},
and further developed in~\cite{FP, post}. 

Let $\E_n$ be the associative algebra 
generated by the symbols $[ij]$, for all $i,j\in\{1,\dots,n\}$, $i\neq
j$, subject to the convention $[ij]+[ji]=0$ and the relations 
\begin{eqnarray}
\label{eq:quadratic}
\begin{array}{l}
[ij]^2 = 0\ ,                    \\[.1in] 
[ij][jk]+[jk][ki]+[ki][ij] = 0\ , \ \ 
\textrm{$i,j,k$ distinct,} \\[.1in] 
[ij][kl]-[kl][ij] = 0\ ,\ \ \textrm{$i,j,k,l$ distinct.} 
\end{array}
\end{eqnarray}
The algebras $\E_n$ are naturally graded; 
the formulas for their Hilbert polynomials, 
for $n\leq 5$, can be found in~\cite{FK}.  
The algebras $\E_n$ are not Koszul for $n\geq 3$ (proved by
Roos~\cite{roos}). 
It is unknown whether $\E_n$ is generally finite-dimensional;
it was proved in~\cite{FP} that the Hilbert series of $\E_n$ divides
that of~$\E_{n+1}\,$.

The ``Dunkl elements'' $\dunkl_1,\dots,\dunkl_n\in\E_n$ are defined by 
\begin{equation}
\label{eq:dunkl1}
\dunkl_j=-\sum_{j<j}[ij]+\sum_{j<k}[jk] \,. 
\end{equation}

\begin{theorem}
\label{th:h*}
\cite{FK}
The complete list of relations satisfied by the Dunkl elements
$\dunkl_1,\dots,\dunkl_n\in\E_n$
is given by $\dunkl_i \dunkl_j=\dunkl_j\dunkl_i$ (for any $i$ and~$j$) and 
$e_i(\dunkl_1,\dots,\dunkl_n)=0$ (for $i=1,\dots,n$).
Thus these elements generate a commutative subring 
canonically isomorphic to $P_n/I_n\,$, and to the cohomology ring of the flag manifold.
\end{theorem}

Let $s_{ij}\in S_n$ 
denote the transposition of elements~$i$ and~$j$.
Consider the ``Bruhat operators''~$[ij]$ 
acting in the group algebra of~$S_n$ by
\begin{equation}
\label{eq:bruhat-operators}
[ij]\,w = 
\left\{
\begin{array}{ll}
w s_{ij} & \textrm{if $\l(w s_{ij})=\l(w)+1$\ ;} \\[.1in]
0           & \textrm{otherwise\ .} \\
\end{array}
\right.
\end{equation}
One easily checks that these operators satisfy
the relations~(\ref{eq:quadratic}).
We thus obtain an (unfaithful) representation of the algebra~$\E_n\,$,
called the \emph{Bruhat representation}.
This representation 
has an equivalent description in the language
of Schubert polynomials.
Let us identify each element $w\in S_n$ with the
corresponding Schubert polynomial~$\S_w\,$.
Then the generators of~$\E_n$ act 
in $\mathbb{Z}[x_1,\dots,x_n]/I_n$ by
\begin{equation}
\label{eq:bruhat-schubert}
[ij]\,\S_w = 
\left\{
\begin{array}{ll}
\S_{w s_{ij}} & \textrm{if $\l(w s_{ij})=\l(w)+1$\ ;} \\[.1in]
0           & \textrm{otherwise\ .} \\
\end{array}
\right.
\end{equation}
The following result is a restatement of the classical
Monk's rule~\cite{monk}.

\begin{theorem}
\label{th:monk}
In the representation {\rm (\ref{eq:bruhat-schubert})} 
of the quadratic algebra~$\E_n$
in the quotient ring $\mathbb{Z}[x_1,\dots,x_n]/I_n\,$, 
a Dunkl element~$\dunkl_j$
acts as multiplication by~$x_j\,$, for $j=1,\dots,n$. 
In other words, $x_j f=\dunkl_j f$, 
for any coset $f\in\mathbb{Z}[x_1,\dots,x_n]/I_n\,$.
\end{theorem}

\section{Structure constants and nonnegativity conjecture}
\label{sec:nonnegativity}

Let $c^w_{uv}$ denote the coefficient of~$\S_w$
in the product~$\S_u\S_v\,$.
Equivalently, $c^w_{uv}$ is 
the number of points in the intersection of the general
translates of three (dual) Schubert cells labelled by~$u$, $v$,
and~$\wnot w$, respectively.
Thus all the~$c^w_{uv}$ are nonnegative integers. 
The problem of finding a combinatorial interpretation for~$c^w_{uv}$
is  one of the central open problems in Schubert calculus. 
In fact, no elementary proof
of the fact that~$c^w_{uv}\ge 0$ is known. 
Much less is known about the more general 
Gromov-Witten invariants of the flag manifold.

Let $\Enp\subset\E_n$ be the cone of all elements that can be written
as nonnegative integer 
combinations of noncommutative monomials
in the generators~$[ij]$, for~$i<j$.

\begin{conjecture}
\label{conj:nonnegativity}
\cite{FK}
\emph{(Nonnegativity conjecture)}
For any $w\in S_n\,$,
the Schubert polynomial~$\S_w$ evaluated at the Dunkl
elements belongs to the positive cone~$\Enp\,$:
\begin{equation}
\label{eq:nonnegativity}
\S_w(\dunkl)=\S_w(\dunkl_1,\dots,\dunkl_{n-1}) \in \Enp\ .
\end{equation}
\end{conjecture}

\setlength{\unitlength}{0.4mm} \thicklines
\begin{figure}
\begin{center}
 \begin{picture}(40,85)(0,-5)
    \put(20,0){\circle*{3}} \put(0,40){\circle*{3}}
    \put(0,20){\circle*{3}} \put(20,60){\circle*{3}}
    \put(40,20){\circle*{3}} \put(40,40){\circle*{3}}
    \put(20.4,-1){\line(1,1){20}}
    \put(19.6,1){\line(1,1){20}}
    \put(20,0){\line(-1,1){20}}
    \put(-1,20){\line(0,1){20}}
    \put(1,20){\line(0,1){20}}
    \put(40,20){\line(0,1){20}}
    \put(0,40){\line(1,1){20}}
    \put(39.4,39.4){\line(-1,1){20}}
    \put(40.6,40.6){\line(-1,1){20}}
    \put(17,-12){$1$}
    \put(-51,16){$[12]+[13]$}
    \put(46,16){$[13]+[23]$}
    \put(-78,39){$[13][23]+[23][13]$}
    \put(46,39){$[12][13]+[13][12]$}
    \put(-54,67){$[12][13][23]+[13][12][13]+[13][23][13]$}
  \end{picture}
\end{center}
  \caption{Evaluations of Schubert polynomials 
at Dunkl elements}
  \label{fig:Schub3(a)}
\end{figure}
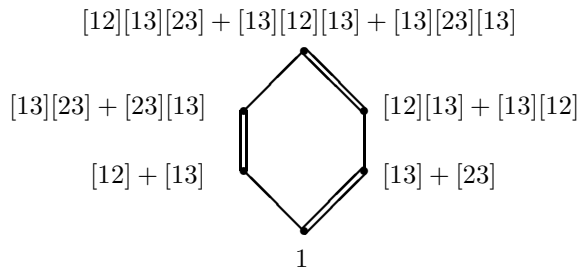

Let us now explain why Conjecture~\ref{conj:nonnegativity}
implies nonnegativity of the structure constants~$c^w_{uv}\,$,
and why furthermore a combinatorial description for the
evaluations~$\S_w(\dunkl)$ 
would provide a combinatorial rule describing the~$c^w_{uv}\,$.

The action~(\ref{eq:bruhat-schubert})
of~$\E_n$ on the quotient ring
$\mathbb{Z}[x_1,\dots,x_n]/I_n$ is defined in such a way
that every noncommutative monomial in the generators
$[ij]$, $i<j$, when applied to a Schubert polynomial~$\S_v\,$,
gives either another Schubert polynomial or zero.
It follows that, for any $z\in\Enp\,$,
the polynomial $z\S_v$ is \emph{Schubert-positive}, i.e.,
is a nonnegative linear combination of Schubert polynomials.
In particular, if Conjecture~\ref{conj:nonnegativity}
holds, then the polynomial $\S_u(\dunkl)\S_v(x)$
is Schubert-positive (here $x$ stands for $x_1,\dots,x_n\,$).
Since, according to Theorem~\ref{th:monk},
\begin{equation}
\label{eq:SaSx=SxSx}
\S_u(\dunkl)\S_v(x)=\S_u(x)\S_v(x)\ ,
\end{equation}
we conclude that $\S_u\S_v$ is Schubert-positive,
i.e., the structure constants $c^w_{uv}$ are nonnegative. 
Now suppose we have a combinatorial description
for~$\S_u(\dunkl)$. By~(\ref{eq:SaSx=SxSx}),
\begin{equation}
\label{eq:cuvw}
c^w_{uv}
= \langle\,\textrm{coefficient of } w
\textrm{ in } \S_u(\dunkl)\,v\,\rangle\ ,
\end{equation}
where the action of $\S_u(\dunkl)$ on $v\in S_n$ is the Bruhat
representation action~(\ref{eq:bruhat-operators}).
Thus~(\ref{eq:cuvw}) would provide a combinatorial rule
for~$c^w_{uv}\,$.


The following conjecture, if proved, would provide an alternative
description of the basis of Schubert cycles. 

\begin{conjecture} \cite{FK}
The evaluations $\S_w(\dunkl)$ are the additive generators of 
the intersection of the cone~$\Enp$
with the commutative subalgebra 
generated by the Dunkl elements.
\end{conjecture}

\section{Quantum deformation of the quadratic algebra}
\label{sec:Eq}

The \emph{quantum deformation}~$\Eqn$ 
of the quadratic algebra~$\E_n$ is defined by replacing the relation
$[ij]^2=0$ in~(\ref{eq:quadratic}) by 
\begin{equation} 
\label{eq:q-i}
[ij]^2=
\left\{
\begin{array}{ll}
q_i & \textrm{if $j=i+1$\ ;} \\
0   & \textrm{otherwise\ .} \\
\end{array}
\right. 
\end{equation} 
The ``quantum Bruhat operators''~$[ij]$, acting in 
the $\mathbb{Z}[q_1,\dots,q_{n-1}]$-span 
of the symmetric group~$S_n$ by
\begin{equation}
\label{eq:q-bruhat-operators}
[ij]\,w = 
\left\{
\begin{array}{ll}
      w s_{ij} & \textrm{if $\l(w s_{ij})=\l(w)+1$\ ;} \\
q_{ij}w s_{ij} & \textrm{if $\l(w s_{ij})=\l(w)-\l(s_{ij})$\ ;} \\
0           & \textrm{otherwise\ ,} \\
\end{array}
\right. 
\end{equation}
provide a representation of $\Eqn\,$, which degenerates into the
ordinary Bruhat representation in the classical limit. 
The operators~{\rm(\ref{eq:q-bruhat-operators})} can be
viewed as acting in the quotient space 
$\mathbb{Z}[q_1,\dots,q_{n-1}][x_1,\dots,x_n]/I^q_n$ by
\begin{equation}
\label{eq:q-bruhat-schubert}
[ij]\,\S_w^q = 
\left\{
\begin{array}{ll}
      \S^q_{w s_{ij}} & \textrm{if $\l(w s_{ij})=\l(w)+1$\ ;} \\
q_{ij}\S^q_{w s_{ij}} & \textrm{if $\l(w s_{ij})=\l(w)-\l(s_{ij})$\ ;} \\
0           & \textrm{otherwise\ .} \\ 
\end{array}
\right. 
\end{equation}

The Dunkl elements~$\dunkl_j\in\Eqn$ are defined by the same 
formula~(\ref{eq:dunkl1}) as before.

\begin{theorem}
\label{th:q-monk}{\rm\cite{FGP}}
{\rm (Quantum Monk's formula)}
In the representation {\rm (\ref{eq:q-bruhat-schubert})} 
of~$\,\Eqn\,$,
a~Dunkl element~$\dunkl_j$
acts as multiplication by~$x_j\,$, for $j=1,\dots,n$.
\end{theorem}

The following result is a corollary of Theorem~\ref{th:q-monk}.

\begin{corollary}
\label{cor:QSP-description}
As an element of the quotient ring~$P_n/I^q_n\,$, 
a quantum Schubert polynomial~$\S^q_w$ 
is uniquely defined by the condition that, in the quantum Bruhat 
representation~{\rm(\ref{eq:q-bruhat-operators})}, it acts on the
identity permutation~$1$ by 
$w = \S^q_w(\dunkl_1,\dots,\dunkl_n)(1)$. 
\end{corollary}

The quantum analogue of Theorem~\ref{th:h*} stated below was conjectured
in~\cite{FK} and proved by A.~Postnikov in~\cite{post}. 

\begin{theorem}
The commutative subring generated by the Dunkl elements
in the quadratic algebra~$\Eqn$
is canonically isomorphic to the quantum cohomology ring
of the flag manifold.
The isomorphism is defined by
$\dunkl_1+\cdots+\dunkl_j \longmapsto \s_{s_j}\,$.
\end{theorem}

The following statement strengthens and refines 
Conjecture~\ref{conj:nonnegativity}.

\begin{conjecture}
\label{conj:q-nonnegativity}
\cite{FK}
For any $w\in S_n\,$, 
the evaluation~$\S^q_w(\dunkl_1,\dots,\dunkl_n)$
can be written as a linear combination of monomials
in the generators~$[ij]$,
with nonnegative integer coefficients.
\end{conjecture}

It is not even clear a priori
that the evaluations~$\S^q_w(\dunkl)$
can be expressed as linear combinations of monomials
with coefficients not depending on the quantum 
parameters $q_1,\dots,q_{n-1}\,$.

A reformulation of (\ref{eq:2.3.1}) in the language
of quantum Schubert polynomials gives 
\begin{equation}
  \label{eq:mult-schubert}
  \S^q_u \S^q_v\,=\,\sum_{w\in S_n} \sum_{d}q^d
  \left<\s_u,\s_v,\s_w\right>_d \S^q_{\wnot w} \ . 
\end{equation}
In view of Theorem~\ref{th:q-monk},
one obtains the following analogue of~(\ref{eq:cuvw}).

\begin{corollary}
\cite{FK}
For $u,v,w\in S_n$ and 
$d=(d_1,\dots,d_{n-1})\in\mathbb{Z}_+^{n-1}\,$,
we have 
\begin{equation}
\label{eq:sss}
\left<\s_u,\s_v,\s_w\right>_d
= \langle\,\textrm{coefficient of } \,q^d\wnot w\,
\textrm{ in } \S^q_u(\dunkl)\,v\,\rangle\ ,
\end{equation}
where $\S^q_u(\dunkl)$ acts on $v$ according to 
the quantum Bruhat representation~{\rm(\ref{eq:q-bruhat-operators})}.
\end{corollary}

Assuming Conjecture~\ref{conj:q-nonnegativity} holds, 
one would like to have a combinatorial rule for 
a nonnegative expansion of~$\S^q_w(\dunkl)$.
Such a rule would immediately lead to a direct combinatorial description
of the Gromov-Witten invariants $\left<\s_u,\s_v,\s_w\right>_d$  of
the flag manifold, given by~(\ref{eq:sss}).

\end{document}